\documentclass{amsart}
\usepackage{amscd}
\usepackage{amssymb}
\usepackage[matrix,knot,all]{xy}
\newcommand{\K}{\mathcal K}
\newcommand{\N}{\mathbb N}
\newcommand{\R}{\mathbb R}
\newcommand{\Q}{\mathbb Q}
\newcommand{\Z}{\mathbb Z}
\newcommand{\tD}{\tilde{\mathcal D}}

\newcommand{\x}{\times}
\newcommand{\map}{Map}

\newtheorem{thm}{Theorem}
\newtheorem{defi}[thm]{Definition}  
\newtheorem{prop}[thm]{Proposition}
\newtheorem{cor}[thm]{Corollary}
\newtheorem{lem}[thm]{Lemma}

\newcommand{\sd}{\rtimes}
\newcommand{\de}{\partial}
\newcommand{\htot}{\widetilde{Tot}}
\newcommand{\hotot}{\widetilde{Tot}}
\newcommand{\ot}{\otimes}
\newcommand{\D}{\mathcal{D}}

\title{Knots, operads and double loop spaces}

\author{Paolo Salvatore}

\address{Dipartimento di matematica, Universit\`a di Roma ``Tor Vergata'', Roma, Italy}
\email{salvator@mat.uniroma2.it}

\begin{document}
\maketitle
\begin{abstract}
We show that the space of long knots in an euclidean space of dimension
larger than three  is a double loop space, proving a conjecture
by Sinha. 
We also construct a double loop space structure on framed long knots,
and show that the map forgetting the framing 
is not a double loop map in odd dimension.
However there is always such a map in the reverse direction expressing the double loop space of framed long knots as a semidirect product.
A similar compatible decomposition  holds for the 
homotopy fiber of the inclusion of long knots into immersions.
We also show via string topology that the space of closed knots in a sphere, suitably desuspended, admits an action of the little 2-discs operad in the category of spectra.
A fundamental tool is the McClure-Smith cosimplicial machinery, that produces double loop spaces out
of topological operads with multiplication.
\end{abstract}

\section{Introduction}

The space $Emb_n$ of long knots in $\R^n$ is the space of embeddings $\R \to \R^n$ that agree with a fixed inclusion of 
a line near infinity. 
The space $Emb_n$ is equipped with the Whitney topology, and it can be identified 
up to homotopy with the subspace of based knots in $S^n$ with fixed derivative at the base point.
The proof that $Emb_2$ is contractible goes back to Smale.  
The components of $Emb_3$ correspond to classical knots. 
The homotopy type of those components has been completely described by Ryan Budney \cite{Bu2}.
For $n>3$ the space $Emb_n$ is connected by Whitney's theorem. The rational homology of $Emb_n$ for $n>3$ has been recently computed by Lambrechts, Turchin and Volic \cite{LTV}.

Rescaling and concatenation defines a natural product on the space
of long knots that is associative up to higher homotopies.
Thus $Emb_n$ is an $A_\infty$-space and in the case $n>3$, being connected,
has the homotopy type of a loop space. The product is homotopy commutative, essentially by 
passing one knot through the other. This suggested that $Emb_n$ could be
(up to weak equivalence) a double loop space.
Budney and Sinha proved that two spaces closely related to $Emb_n$ are double loop spaces, for $n>3$, by  different approaches.
A framed long knot in $\R^n$ is a long knot in $\R^n$ together with a choice of framing $\R \to SO(n)$, standard near infinity,
 such that the first
vector of the framing gives the unit tangent vector map $\R \to S^{n-1}$ of the knot.
Budney shows in \cite{Bu} that the space $fEmb_n$
of framed long knots in $\R^n$ is  a double loop space for $n>3$.
This is achieved by constructing an explicit action of the little 2-cubes operad on a space homotopy equivalent to 
the group-like space $fEmb_n$.
The operad action is also defined for $n=3$, and makes $fEmb_3$ into a free 2-cubes algebra on the non-connected space of prime long knots.

Sinha shows in \cite{sinha} that the homotopy fiber $Emb'_n$ of the unit tangent vector map
 $Emb_n \to \Omega S^{n-1}$ is a double loop space, and the map 
is nullhomotopic. His approach goes via the cosimplicial machinery by McClure and Smith \cite{MS}
that produces double loop spaces out of non-symmetric operads in based spaces.
Under this correspondence $Emb'_n$ is produced by an operad equivalent to the little $n$-discs operad, the Kontsevich operad.
We show 
that this machinery, applied to an operad equivalent to the framed little $n$-discs operad,
gives a double loop space structure on framed long knots in $\R^n$,
that presumably coincides with the one described by Budney.
We believe that the fact that the framed little discs is a cyclic operad 
\cite{Bu3} together with the McClure-Smith machinery for cyclic objects 
will lead to
a {\it framed} little 2-discs action on framed long knots. 

Let us consider the principal fibration $$\Omega SO(n-1) \to fEmb_n \to Emb_n$$
forgetting the framing.
Such fibration is trivial because its classifying map $Emb_n \to SO(n-1)$ is the composite of the 
(nullhomotopic) unit tangent vector map and
the holonomy 
$\Omega S^{n-1} \to SO(n-1)$.
Given the splittings  
\begin{equation} \label{'}
Emb'_n \simeq Emb_n \times \Omega^2 S^{n-1}
\end{equation}
 and
\begin{equation} \label{f}
fEmb_n \simeq Emb_n \times \Omega SO(n-1) \end{equation}
 Sinha asked in \cite{sinha} whether one could restrict the double loop structure to the first factor. We answer this affirmatively.

\begin{thm} \label{main}
The space $Emb_n$ of long knots in $\R^n$ is a double loop space for $n>3$. 
\end{thm}

The double loop space structure is not produced directly from an operad as hoped in \cite{sinha}, but is deduced by diagram chasing on a diagram of cosimplicial spaces.

The splittings (\ref{'}) and (\ref{f}) respect the single loop space structures but not the double loop space structures, as the projections on the 
factor $Emb_n$ are not double loop maps in general.

\begin{thm} \label{not2}
The map forgetting the framing $fEmb_n \to Emb_n$
and the map from the homotopy fiber $Emb'_n \to Emb_n$  are not double loop maps for $n$ odd.
\end{thm}

We prove this by showing that the maps in question do not preserve the Browder operation, a natural bracket on the homology of double loop spaces. This is based on computations by Turchin \cite{T}.

There are instead double loop maps
$Emb_n \to Emb'_n$ and $Emb_n \to fEmb_n$ that 
together with the fiber inclusions 
$\Omega^2 S^{n-1} \to Emb'_n$ and $\Omega SO(n-1) \to fEmb_n$
produce
 essentially semidirect product extensions of double loop spaces. We state this precisely
in the following theorem.

\begin{thm}        \label{frame}

There is a commutative diagram of double loop spaces and double loop maps
$$\xymatrix{
Emb_n \ar[r] \ar@{=}[d] &  Emb'_n \ar@<1ex>[r] \ar[d] & \Omega^2 S^{n-1} \ar[d] \ar@<1ex>[l]  \\
Emb_n \ar[r]     &  fEmb_n \ar@<1ex>[r] & \Omega SO(n-1) \ar@<1ex>[l]  }$$

The rows deloop twice to fibrations with sections, and the vertical maps
are induced by the holonomy $\Omega S^{n-1} \to SO(n-1)$.

\end{thm}

Also this theorem develops the approach by Sinha.
The double loop spaces and double loop maps are produced by applying 
the McClure-Smith machinery to 
suitable operads and operad maps. 
  
At the end of the paper we apply ideas from string topology to
show that the shifted homology
of the space $Emb(S^1,S^n)$ of all knots in the $n$-sphere behaves as
the homology of a double loop space. More precisely this structure is induced by the action of an operad equivalent to the little 2-cubes
at the spectrum level rather than at the space level.

\begin{thm} \label{sphere}
The spectrum $\Sigma^{1-2n} \Sigma^\infty Emb(S^1,S^n)_+$ is an $E_2$-ring spectrum.
\end{thm}

A similar result has been obtained independently by Abbaspour-Chataur-Kallel.
The case $n=3$ is joint work with Kate Gruher \cite{GS}.

\medskip

Here is a plan of the paper: 
in section \ref{two} we recall some background material on operads, 
cosimplicial spaces and prove theorem \ref{main}.
In section \ref{three} we study the space of framed knots via cosimplicial
techniques and prove theorem \ref{frame}.
In section \ref{four} we recollect some material on the Deligne conjecture  and we give a proof of theorem \ref{not2}.
In the last section \ref{five} we develop the string topology of knots 
proving theorem \ref{sphere}.

\medskip

I would like to thank Ryan Budney, Pascal Lambrechts, Riccardo Longoni, Dev Sinha and Victor Turchin for helpful conversations regarding this material.

\section{Cosimplicial spaces and knots} \label{two}

We recall that a topological operad $O$ is a collection of spaces
$O(k), \,k \geq 0$, 
 together with
a unit $\iota \in  O(1)$ and composition maps
$$\circ_t:O(k)\times O(l) \to O(k+l-1)\; $$
for $1\leq t \leq k$
 satisfying appropriate axioms \cite{may2}.
The operad is {\em symmetric} if the symmetric group  
$\Sigma_k$ acts on $O(k)$ for each $k$, compatibly with the composition maps.
We say that a space $A$ is acted on by an operad $O$, or it is a
$O$-algebra, if we are given maps
$O(n) \times A^n \to A$ satisfying appropriate associativity and unit 
axioms \cite{may2}.
The concepts of (symmetric) operads and their algebras can be defined likewise in any (symmetric) monoidal caegory.

\

Let $F(\R^n,k)$ be the ordered configuration space of $k$ points  in $\R^n$. 
The direction maps $\theta_{ij}:F(\R^n,k) \to S^{n-1}$  are defined for $i \neq j$ by
$$\theta_{ij}(x_1,\dots,x_n)=(x_i-x_j)/|x_i-x_j|.$$ 
Let us write $B_n(k) = (S^{n-1})^{k(k-1)/2}$. We can think of $B_n(k)$ as the space of formal 
'directions' between $k$ distinct points in
$\R^n$, 
where the directions
are indexed by distinct pairs of integers between 1 and $k$. 
By convention we set $B_n(1)$ and $B_n(0)$ equal to a point.
\begin{prop} \cite{sinha}
The collection $B_n(k)$ forms a symmetric topological operad.  
\end{prop}
The action of the symmetric group $\Sigma_k$ on $B_n(k)$ permutes both indices.
Intuitively the operad composition replaces a point by an infinitesimal configuration and relabels the points.
More precisely we must specify the composition rule $\circ_t:B_n(k) \times B_n(l) \to B_n(k+l-1)$ for $1 \leq t \leq k$.
For elements $\alpha=(\alpha_{ij})_{1\leq i<j\leq k} $  and $\beta=( \beta_{ij})_{1\leq i<j\leq l}$ the composition is
$$(\alpha \circ_t  \beta)_{ij}=   
\begin{cases}
\alpha_{ij} \; {\rm for}\; i<j\leq t    \\
\beta_{i-t+1,j-t+1} \;{\rm for}\; t\leq i < j \leq t+l-1 \\
\alpha_{i-l+1,j-l+1}\;   {\rm for}\; t+l \leq i<j \\
\alpha_{i,t} \;  {\rm for}\; i<t \leq j < t+l  \\
\alpha_{t,j} \;  {\rm for}\; t \leq i < t+l \leq j 
\end{cases}.$$

Let $\theta^k: F(\R^n,k) \to B_n(k)$ be the product of all direction maps 
$\theta^k(x)=(\theta_{ij}(x))_{1\leq i <j \leq n}.$
For $k\geq 2$ let $\K_n(k) \subset B_n(k)$ be the closure of the image of $\theta^k$.
We set also $\K_n(0)=B_n(0)=\{*\}$ and $\K_n(1)=B_n(1)=\{\iota\}$.
The restriction $\theta^k:F(\R^n,k) \to \K_n(k)$ is a $\Sigma_k$-equivariant homotopy equivalence.

\begin{prop} \cite{sinha}
The collection $\K_n(k)$ forms a suboperad of $B_n(k)$ that is weakly equivalent to the little $n$-discs operad. \end{prop}
The operad $\K_n$ is known as the Kontsevich operad.

We say that a non-symmetric topological operad has a multiplication if there is a choice of base points $m_k \in O(k)$ for each $k$ such that the structure
maps are based maps. This is the same as a non-symmetric
operad in based spaces.

The operads $B_n$ and $\K_n$ have a multiplication, 
defined by setting all components $\theta_{ij}\, (i<j)$ of the base points $m_k$ equal to a fixed direction. We choose the
last vector of the canonical basis of $\R^n$ as fixed direction.

We recall the definition of a cosimplicial space.
Let $\Delta$ be the category with standard ordered sets $[k]=\{0<\dots<k\}$ as objects $(k \in \N)$ 
 and monotone maps as morphisms.
A cosimplicial space is a covariant functor from the category $\Delta$ to the category of topological spaces.
For each $k$ the simplicial set $\Delta(\_,[k])$ is also called the simplicial $k$-simplex $\Delta^k_*$ \, .
Its geometric realization is the standard $k$-simplex $\Delta^k$. All simplexes fit together to form a cosimplicial space.
In fact if we apply geometric realization to the bisimplicial set (functor from $\Delta$ to simplicial sets) 
  $\Delta(*',*)$ in the variable $*'$ then we obtain a cosimplicial space denoted by
$\Delta^*$.

The totalization $Tot(S^*)$ of a cosimplicial space $S^*$ is the space of natural transformations
 $\Delta^* \to S^*$. There is a standard cosimplicial map $\tilde{\Delta^*} \to \Delta^*$, where $\tilde{\Delta}^*$ is an appropriate cofibrant resolution.
The
homotopy totalization $\htot(S^*)$ is the space of natural transformations
$\tilde{\Delta}^* \to S^*$.  
This is also the homotopy limit of the functor from $\Delta$ to spaces defining the cosimplicial space.
Precomposition induces a canonical map $Tot(S^*) \to \htot(S^*)$ that is a weak equivalence when $S^*$ is fibrant,
in the sense that it satisfies the matching condition \cite{hirschhorn}.

An operad $(O,p)$ with multiplication defines a cosimplicial space $O^*$ sending $[k]$ to $O(k)$.
The cofaces operator $d^i:O(k) \to O(k+1)$ is defined by  
$$\begin{cases}
d^i(x)=x \circ_i m_2\;{\rm for}\; 1\leq i \leq k  \\
d^0(x)=m_2 \circ_1 x \\
d^{n+1}(x)=m_2 \circ_2 x.
\end{cases} $$ 
The codegeneracies $s^i:O(k)\to O(k-1)$  are defined by
$s^i(x)=x \circ_i m_0$.

\begin{thm} {\rm (McClure-Smith)} \label{ms}
Let $O$ be an operad with multiplication.
Then the totalization $Tot(O^*)$ 
(respectively the homotopy totalization $\widetilde{Tot}(O^*)\;$) admits
 an action of an operad $\D_2$ (respectively $\tD_2$) weakly equivalent to the little 2-cubes operad.
\end{thm} 
 
By the recognition principle \cite{may2}
if $Tot(O^*)$ or $\widetilde{Tot}(O^*)$
is connected then it is weakly equivalent to a double loop space.

\

Given a simplicial set $S_*$, considered as simplicial space with discrete values,
and a space $X$, we obtain a cosimplicial space
$map(S_*,X)$, often denoted $X^{S_*}$. If $S_*$ is a simplicial based set and $X$ is a based space 
then we obtain similarly a cosimplicial space
$map_\bullet(X_*,S)$. 
Let us denote by $|S|$ the geometric realization of $S$.
The following is standard.

\begin{prop} \label{homeo}
The adjoint maps of the evaluation maps 
$$map(|S|,X) \times \Delta_k \to map(S_k,X)$$
induce a homeomorphism $map(|S|,X) \cong Tot(map(S_*,X))$.
In the based version we obtain a homeomorphism from the based mapping space
$$map_\bullet(|S|,X) \to Tot(map_\bullet(S_*,X)).$$
 The canonical maps from these totalizations to the homotopy 
totalizations are weak equivalences.

\end{prop}

Let $\Delta^k_*$ be the simplicial $k$-simplex, and $\de \Delta^k_*$ its simplicial subset obtained by removing
the non-degenerate simplex in dimension $k$ and its degeneracies. 
The quotient $S^k_*:=\Delta^k_*/\de \Delta^k_*$ is the simplicial $k$-sphere.

\begin{prop} \cite{sinha} \label{bn}

The cosimplicial space $B_n^*$ is isomorphic to $\map_\bullet(S^2_*,S^{n-1}).$
\end{prop}

Namely $B_n^k$ has a factor $S^{n-1}$ for each pair $1\leq i < j \leq k$ and
$map_\bullet(S^2_*,S^{n-1})$ has a sphere factor for each $k$-simplex of $S^2_*$, namely 
for each non-decreasing sequence 
of length $k+1$ starting with 0 and ending with 2. Then $i$ corresponds to the position of the last $0$ 
and $j$ to the position of the last $1$.
Propositions \ref{homeo} and \ref{bn} imply the following corollary.

\begin{cor}
The totalization $Tot(B_n^*)$ is homeomorphic to $\Omega^2(S^{n-1})$.
\end{cor}

There is also a cosimplicial space $\K_n^* \sd S^{n-1}$, not defined by an operad with multiplication.
This is constructed so that $\K_n^k \sd S^{n-1} =\K_n(k) \times (S^{n-1})^k $.
Elements of this space can be thought of as configurations of $k$ points in $\R^n$, each labelled by a  direction.
The composition rule can be defined as follows, via the identification $S^{n-1} = \K_n(2)$.
Given $(x;v_1,\dots,v_k) \in \K_n(k) \times (S^{n-1})^k$,  we define
for $1 \leq i \leq k$
$$d^i(x;v_1,\dots,v_k)=(x\circ_i v_i;v_1,\dots,v_i,v_i,\dots,v_k).$$
Intuitively  these cofaces double a point in the associated direction, at infinitesimal distance.
The first and last cofaces add a point labelled by the preferred direction 'before' or 'after' 
the
configuration and are defined by
 $$d^0(x;v_1,\dots,v_k)=(m_2 \circ_1 x;v_1,\dots,v_k,m_2)$$ and
$$d^{k+1}(x;v_1,\dots,v_k)=(m_2 \circ_2 x;m_2,v_1,\dots,v_k).$$
The codegeneracies forget a point and are defined by
$$s^i(x;v_1,\dots,v_k)=(x\circ_i m_0;v_1,\dots,\hat{v_i},\dots,v_k).$$
The very same rules define a cosimplicial space
 $B_n^* \sd S^{n-1}$ with
$B_n^k \sd S^{n-1} = (S^{n-1})^{k(k-1)/2} \times (S^{n-1})^k$ so that
$\K_n^* \sd S^{n-1} \subset B_n^* \sd S^{n-1}$ is a cosimplicial subspace.

\begin{thm} {\rm (Sinha)} \cite{sinha} \label{embn}
The homotopy totalization of 
 $\K_n^* \sd S^{n-1}$ is weakly equivalent to
$Emb_n$.
\end{thm}
The proof of this theorem relies on Goodwillie calculus.
From now on we will mean by $Emb_n$ the 
space of smooth maps from the interval
$I$ to
the cube $I^n$ sending the extreme points of the interval to centers of opposite faces of the cube, 
with derivative orthogonal to the faces.

The weak equivalence $Emb_n \to \htot(\K_n^* \sd S^{n-1})$ is constructed as follows, by evaluating directions between points of the knot
and tangents.
Regard an element of the $k$-simplex as a sequence of real numbers $0 \leq x_1 \leq \dots \leq x_k  \leq 1$.
There are maps $\beta_k:Emb_n \to map(\Delta^k,\K_n(k) \times (S^{n-1})^k)$ defined by
$$\beta_k(f)(x_1,\dots,x_k)=\{ \theta^k(f(x_1),\dots,f(x_k)),
 f'(x_1)/|f'(x_1),\dots, \,f'(x_k)/|f'(x_k) \}$$
when $x_1 <\dots <x_k$. If some $x_i=x_j$ for $i<j$ then we must replace the component $\theta_{ij}=f(x_j)-f(x_i)/|f(x_j)-f(x_i)| $ in the expression above by 
 $f'(x_i)/|f'(x_i)|$. All maps $\beta_k$ fit together to define a map
$\beta:Emb_n \to Tot(\K_n^* \sd S^{n-1})$. The composite with the standard map to the homotopy totalization is the desired 
weak equivalence.

\

Let us recall some background on homotopy fibers:
the homotopy fiber of a based map $f:X \to Y$ is defined by the pullback square
$$\begin{CD}   Hofib(f)             @>>>       X \\
              @VVV                                          @VVfV               \\ 
                PY            @>ev>>            Y
\end{CD}$$
with $PY$ the contractible space of paths in $Y$ sending 0 to the base point, and $ev$ the evaluation at the point 1.
If $f$ is a fibration with fiber $F$ then there is a canonical homotopy equivalence $F \to Hofib(f)$ sending $x \in F \subset X$ to the 
pair $(x,c)$ with $c$ the constant loop at the base point of $Y$. 
The homotopy fiber is homotopy invariant,
namely given a commutative diagram 
$$\begin{CD}   X             @>f>>       Y  \\
              @VV\simeq V                                          @V\simeq VV               \\ 
               X'            @>f'>>      Y'
\end{CD}$$
with the vertical arrows weak equivalences, then the induced map
$Hofib(f) \to Hofib(f')$ is a weak equivalence.
This is a special case of the homotopy invariance of homotopy limits ( theorem 18.5.3 (2) in \cite{hirschhorn}).

\
\begin{cor} {\rm (Sinha)}
The homotopy fiber $Emb'_n$ of the unit tangent vector map $u: Emb_n \to \Omega S^{n-1}$ is weakly equivalent to the homotopy totalization $\widetilde{Tot}(K_n^*)$  ,  and
thus is
a double loop space for $n>3$.
\end{cor}

Proof: 
The projection $\K_n(k) \times (S^{n-1})^k \to (S^{n-1})^k$ defines a map of 
cosimplicial spaces
$\K_n^* \sd S^{n-1} \to map_\bullet(S^1_*, S^{n-1})$ and 
there is a commutative square
$$\begin{CD}  Emb_n             @>>\tilde{\beta}>            \htot(\K_n^* \sd S^{n-1}) \\
              @VV{u}V                                          @VV{\pi}V               \\ 
               \Omega S^{n-1}             @>>>            \htot(map_\bullet(S^1_*, S^{n-1})).
\end{CD}$$
By theorem 18.5 (2)
 in \cite{hirschhorn} the homotopy totalization of a sequence of cosimplicial spaces $X^* \to Y^* \to Z^*$ that
are levelwise fibrations is a fibration $\hotot X^*  \to \hotot Y^* \to \hotot Z^*$.
  
Then we have a weak equivalence $\htot(\K_n^*) \to Hofib(\pi)$ and
by homotopy invariance weak equivalences $Emb'_n = Hofib(u) \simeq Hofib(\pi) \simeq  \htot(\K_n^*)$.
We conclude by theorem \ref{ms}. 
 
\

{\em Remark}:  We may substitute $\Omega S^{n-1}$ in the statement above by the space $Imm(I,I^n)$ of 
immersions $I \to I^n$ with fixed values and tangent vectors at the boundary, and $u$ by the inclusion $Emb(I,I^n) \to Imm(I,I^n)$,
 because the unit tangent vector map induces
the Smale homotopy equivalence $Imm(I,I^n) \simeq \Omega S^{n-1}$. 

\

In the next lemma we identify
the totalization of $B_n^* \sd S^{n-1}$.
There are standard simplicial inclusions $d^0_*:\Delta^1_* \to \Delta^2_*$ and 
$d^2_*:\Delta^1_* \to \Delta^2_*$ induced by
 the strictly monotone maps $[1]\to[2]$ avoiding 
respectively 2 and 0.

\begin{lem}
The totalization of the levelwise fibration of cosimplicial spaces
$$B_n^* \to B_n^* \sd S^{n-1} \to    map(\Delta^1_* / \partial \Delta^1_*, S^{n-1})$$
is the fibration 
$$map_\bullet(\Delta^2/\de \Delta^2, S^{n-1}) \to map_\bullet(\Delta^2/(d^0(\Delta^1) \cup d^2(\Delta^1)),S^{n-1}) \to  map_\bullet(\Delta^1/ \de \Delta^1,S^{n-1}).$$

\end{lem}

\begin{proof} 
The space $B_n^k \sd S^{n-1}$ has a factor $S^{n-1}$ for each pair $1\leq i < j \leq k$ and a factor $S^{n-1}$ for each $1 \leq l \leq k$.
The space $$map_\bullet(\Delta^2_k/d^0_k(\Delta^1_k) \cup d^2_k(\Delta^1_k),S^{n-1})$$
has a factor $S^{n-1}$ for each non-decreasing sequence
of length $k+1$ containing 0,1,2
and a factor $S^{n-1}$
 for each non-decreasing sequence of length $k+1$ starting with 0, ending with 2, without 1's.
For these latter sequences
$l$ corresponds to the
last position containing a 0. For the former sequences we apply the same correspondence as 
in the proof of proposition \ref{bn}.
\end{proof}
\

{\em Proof of theorem \ref{main}}.

If we map the sequence $\K_n^* \to \K_n^* \sd S^{n-1} \to map_\bullet(S^1_*,S^{n-1})$ to the sequence
 $B_n^* \to B_n^* \sd S^{n-1} \to map_\bullet(S^1_*,S^{n-1})$ we    
obtain a commutative diagram of cosimplicial spaces that at level $k$ is

$$\begin{CD}  \K_n(k)              @>>>                  (S^{n-1})^{k(k-1)/2} \\
              @VVV                                          @VVV               \\ 
             \K_n(k) \times (S^{n-1})^k  @>>>            (S^{n-1})^{k(k-1)/2} \times (S^{n-1})^k \\ 
             @VVV                                          @VVV                \\
               (S^{n-1})^k         @=                       (S^{n-1})^k .
\end{CD}$$

The homotopy totalization functor gives a diagram of spaces weakly equivalent to those in the 
diagram

$$\begin{CD}   Emb'_n          @>>>                 \Omega^2 S^{n-1}          \\           
                @VVV                                       @VVV                   \\
               Emb_n            @>>>                 P\Omega S^{n-1}            \\ 
                @VVV                                        @VVV                     \\
                \Omega S^{n-1}         @=                  \Omega S^{n-1}.       
\end{CD} $$

Let us analyze the diagram of homotopy totalizations.
By naturality the upper row is a map of algebras over the McClure-Smith operad $\tD_2$, and
then its homotopy fiber $F$ is also an algebra over $\tD_2$.
The homotopy fiber of the second row is weakly equivalent to $Emb_n$ by theorem \ref{embn} and because the target is contractible.  
The homotopy fiber of the third row is contractible. 
The homotopy fibers of the rows in a diagram whose columns are fibrations 
form a fibration ( 18.5.1 in \cite{hirschhorn}),
so that $F \simeq Emb_n$.  This space is connected by Whitney's theorem for $n>3$,
and then
by the recognition principle \cite{may2} 
is weakly equivalent to a double loop space. $\Box$

\section{Framed knots and double loop fibrations} \label{three}

We start by some general considerations on framed knots. By definition $fEmb_n$ is the pullback

$$\begin{CD}   fEmb_n          @>>>                  Emb_n              \\           
                @VVV                                       @V{u}VV                   \\                                                        
                \Omega SO(n)         @>>>                  \Omega S^{n-1}.       
\end{CD} $$
Actually  $fEmb_n$ is homeomorphic to the homotopy fiber of the composite 
$$Emb_n \stackrel{u}{\to} \Omega S^{n-1} \stackrel{h}{\to} SO(n-1)$$
of the holonomy $h$ and the unit tangent vector map $u$.
 The homeomorphism is induced by the projection 
$fEmb_n \to Emb_n$ and the map $fEmb_n \to PSO(n-1)$ considering the difference 
between the framing induced by the holonomy along the knot and the assigned framing of the framed knot.
By naturality of the homotopy fiber construction the holonomy induces a map $Emb'_n \to fEmb_n$.

\

We will give next an operadic interpretation of framed knots.
We recall \cite{SW} that a topological group $G$ acts on a topological operad $O$ if each $O(n)$ is a $G$-space and the operadic
composition maps
are $G$-equivariant. In other words $O$ is an operad in the category of $G$-spaces.
 In such case one can define the semidirect product \cite{Markl,SW}
  $O \sd G$ with $n$-ary space
$O(n) \times G^n$ and composition
$$(p;g_1,\dots,g_n) \circ_i (q;h_1,\dots,h_m) =(p \circ_i g_i(q);g_1,\dots,g_ih_1,\dots,g_ih_m,\dots,g_n).$$

For example, the (trivial) action of a group $G$ on the commutative operad $Com$ 
defines a semidirect product $\underline{G} := Com \sd G$ such that
 $\underline{G}(n)=G^n$.
The framed little $n$-discs operad is isomorphic to the semidirect
$fD_n=D_n \sd SO(n)$, where $SO(n)$ rotates the picture of the little discs.

The natural action of $SO(n)$ on $S^{n-1}$ defines a $SO(n)$-action on the operad $B_n$, given that $B_n(k)=(S^{n-1})^{k(k-1)/2}$.
This action restricts to an action on the operad
$\K_n$. The arguments giving the weak equivalence between $\K_n$ and the little $n$-discs operad $D_n$ extend to show that
the semidirect product operad $f\K_n= \K_n \sd SO(n)$
is weakly equivalent to the framed little $n$-discs operad.
Namely in \cite{barcelona} we constructed a diagram of weak equivalences of operads $D_n \leftarrow WD_n \to F_n$, where $F_n$ is the Fulton-MacPherson 
operad. These arrows and the projection $F_n \to K_n$, that is also a weak equivalence \cite{sinha}  are $SO(n)$-equivariant.

\begin{prop}
The homotopy totalization of the cosimplicial space $f\K_n^*$, for $n>3$,
is  weakly equivalent to the space $fEmb_n$ of framed long knots in $\R^n$.
\end{prop}

\begin{proof}
 
The sequence of cosimplicial spaces
$$map_\bullet(S^1_*,SO(n-1)) \to f\K_n ^* \to \K_n^* \sd S^{n-1}$$
 is levelwise
the fibration
$SO(n-1)^k \to SO(n)^k \times \K_n(k) \to (S^{n-1})^k \times \K_n(k)$. 
There is a commutative diagram
$$\xymatrix{
\Omega SO(n-1) \ar[d] \ar[r] & fEmb_n \ar[d]^{\tilde{f\beta}} \ar[r] &  Emb_n \ar[d]^{\widetilde{\beta}} \\
\widetilde{Tot}(map_\bullet(S^1_*,SO(n-1))) \ar[r] & \widetilde{Tot}(f\K_n ^*) \ar[r] & \widetilde{Tot}(K_n^* \sd S^{n-1}) 
}$$
where the rows are fibrations.
The middle arrow $\widetilde{f\beta}$ is the composite of a map $f\beta$ and the canonical
map 
$Tot(f\K_n ^*) \to \widetilde{Tot}(f\K_n ^*)$, where $f\beta$ is adjoint 
to a collection of maps $fEmb_n \times \Delta_k \to \K_n(k) \times SO(n)^k$ that
evaluate directions between points of the framed knot as before and in addition 
evaluate the framings at those points.

The left and right vertical maps
are weak equivalences, and hence the middle vertical map $\widetilde{f\beta}$ is a weak equivalence.

\end{proof}

Now $f\K_n$ is an operad with multiplication, so that
$\hotot (f\K_n)$ has an action of the McClure-Smith operad $\tD_2$ by theorem \ref{ms}.
The space $fEmb_n \simeq Emb_n \times \Omega SO(n-1)$ is grouplike for $n>3$, in the sense
that its components form a group, namely $\Z_2$.
 By the recognition principle \cite{may2} we readily obtain :
\begin{cor}
The space of framed long knots in $\R^n$  is weakly equivalent to a double loop space for $n>3$.
\end{cor}

This recovers the result by Budney \cite{Bu}.

\

We characterize next the semidirect product operad $B_n \sd SO(n)$, that we will also call $fB_n$.
We observe  that
there is an operad inclusion $i_n:\underline{SO(n-1)} \to fB_n$ that we define next. 
Let us identify $SO(n-1)$ to the subgroup of $SO(n)$ fixing the preferred direction $m_2 \in S^{n-1}=B_n(2)$.
We recall that $m_k \in B_n(k)$ is the base point.
Then $i_n(k)$ sends $(g_1,\dots,g_k)\in SO(n-1)^k$ to
$(m_k,g_1,\dots,g_k) \in B_n(k) \times SO(n)^k$. We visualize the image as a configuration of points on a line parallel to the preferred direction,
with the assigned framings. Clearly $i_n$ factors through the operad $f\K_n$.
We remark that $i_n$ does not extend to a section $\underline{SO(n)} \to fB_n$ of the projection $fB_n \to \underline{SO(n)}$.

\begin{prop} \label{equi}

The map $i_n: \underline{SO(n-1)} \to fB_n$  
induces on the (homotopy) totalizations of the associated cosimplicial spaces
a homotopy equivalence that is a double loop map,
so that 

$$\Omega SO(n-1) \simeq Tot(fB_n^*).$$
\end{prop}

\begin{proof}

We have a pullback diagram of cosimplicial spaces
$$\begin{CD}
fB_n^*  @>>>B_n^* \sd S^{n-1} \\
@VVV      @VVV                \\
map_\bullet(S^1_*,SO(n)) @>>> map_\bullet(S^1_*,S^{n-1}).
\end{CD} $$
On totalizations we obtain the pullback diagram
$$\begin{CD}
Tot(fB_n^*) @>>> P\Omega S^{n-1} \\
@VVV      @VVV                \\
\Omega SO(n) @>>> \Omega S^{n-1}.
\end{CD} $$

The inclusion $\underline{SO(n-1)} \to fB_n$ induces on totalizations the standard
homotopy equivalence  
from $\Omega SO(n-1)$
to $Tot(fB_n^*)$, the homotopy fiber of the looped projection 
$\Omega SO(n) \to \Omega S^{n-1}$.
We can replace totalizations by homotopy totalizations in the proposition
since all cosimplicial spaces involved are fibrant.

\end{proof}

{\em Proof of theorem  \ref{frame}}:
We have a diagram of operads

$$\xymatrix{
                & \K_n \ar[d] \ar[r] & B_n  \ar[d]\\
\underline{SO(n-1)} \ar[r] \ar[dr]   & f\K_n \ar[d] \ar[r]& fB_n \ar[dl] \\
                       &     \underline{SO(n)}&             
}$$

The operad inclusion $f\K_n \to  fB_n$ gives on homotopy totalizations, by naturality of the McClure-Smith construction, 
a map of $\tD_2$-algebras $\hotot(f\K_n) \to \hotot(fB_n)$, 
that by naturality of the recognition principle
is a double loop map. 
Its homotopy fiber $F$ 
 is weakly equivalent to 
$Emb_n$ as double loop space, by comparison with the homotopy fiber of  $\hotot(\K_n) \to \hotot(B_n)$ and by the arguments in the proof
of theorem \ref{main}. The double loop map
$\hotot(f\K_n) \to \hotot(fB_n)$
has a double loop section 
because the operad inclusion $\underline{SO(n-1)} \to f\K_n \to fB_n$ induces a weak equivalence that is a double loop map on homotopy totalizations  (proposition \ref{equi}).
This gives the fiber sequence of double loop maps with section
$$Emb_n \to fEmb_n \stackrel{\rightarrow}{\leftarrow} \Omega SO(n-1).$$

Now there is a commutative diagram

$$\xymatrix{
 \Omega SO(n-1) \ar^{j}[r] \ar^{\simeq}[d] & fEmb_n \ar[r] \ar^{\simeq}[d] & \Omega SO(n) \ar^{\simeq}[d]\\
\hotot(SO(n-1)^*) \ar[r] & \hotot(f\K_n^*)\ar[r] & \hotot(SO(n)^*)    
}$$

and the inclusion $j:\Omega SO(n-1) \subset fEmb_n$ represents the subspace of all framings of the trivial knot.
We conclude the proof
by taking homotopy fibers over
$\hotot(SO(n)^*)$.

Namely the homotopy fiber $K'$ of $\hotot{f\K_n^*} \to \hotot{SO(n)^*}$ 
(resp. $B'$ of  $\hotot{fB_n^*} \to \hotot{SO(n)^*}$ )
is canonically
weakly equivalent to $\hotot{\K_n^*} \simeq  Emb'_n$
(resp. to $Tot(B_n^*) \simeq \Omega^2 S^{n-1}$).
Let $\Omega'$ be the homotopy fiber of
$\hotot SO(n-1)^* \to \hotot SO(n)^*$, canonically weakly eqivalent
to $\Omega^2 S^{n-1}$ as double loop space. Then the double loop map
$K' \to B'$ has a double loop section 
because the composite $\Omega' \to K' \to B'$ is a weak equivalence and a double loop map.
This gives the fiber sequence of double loop maps with section
$$Emb_n \to Emb'_n \stackrel{\rightarrow}{\leftarrow} \Omega^2 S^{n-1}.\quad \Box$$

\section{An obstruction to double loop maps} \label{four}

In this section we will prove theorem \ref{not2} by showing that
the projection $fEmb_n \to Emb_n$ from framed knots to knots
and the map $p:Emb'_n \to Emb_n$ from section \ref{two}
do not preserve the Browder operation in rational homology for $n$ odd.
We need to review some notions on homology operations of double loop spaces.

\begin{defi}
An $n$-algebra is an algebra over the homology operad of the little $n$-discs operad.
\end{defi}

In particular a 2-algebra is called a Gerstenhaber algebra. A (graded) $n$-algebra $A$ for $n>1$ is described by assigning a product and a bracket
$$\_*\_: A_i \otimes A_j \to A_{i+j}$$ 
$$[\_,\_]: A_{i} \ot A_j \to A_{i+j+n-1}$$
that satisfy essentially the axioms of a Poisson algebra,
except for signs. We refer to \cite{SW} 
for a full definition.
The action of the little $n$-discs operad on an $n$-fold loop space gives a natural $n$-algebra structure on its homology, such that
the product is the Pontrjagin product and the bracket is called the Browder operation.
In particular the homologies of the double loop spaces $Emb_n, Emb'_n$
and $fEmb_n$  have a natural structure of Gerstenhaber algebras. 

\
Originally Gerstenhaber introduced the algebraic structure bearing his name while studying the Hochschild complex of associative algebras. 
More generally Gerstenhaber and Voronov introduced 
this structure on the Hochschild homology of 
an operad with multiplication in vector spaces.
Let $O$ be an operad in vector spaces together with a multiplication,
i.e. an operad map $Ass \to O$ from the associative operad.
The image of the multiplication in $Ass$ is an element $m \in O(2)$.
The operad composition maps define a bracket
$$[\_,\_]:O(k) \ot O(l) \to O(k+l-1)$$
by
$$[x,y]=\sum_{i=1}^k \pm  x \circ_i y -  \sum_{i=1}^l  
\pm y \circ_i x$$ 
for appropriate signs \cite{T}.
The multiplication defines a star product 
$$\_ * \_ : O(k) \ot O(l) \to O(k+l)$$
 by
$$x * y = m(x,y) := (m\circ_2 y)\circ_1 x .$$

\begin{defi}
The Hochschild complex of $O$ is the chain complex 
$(\bigoplus s^{-k}O(k), \de)$, where $s^{-k}$ is degree desuspension, and
the differential is $\de(x)=[m,x]$.
The Hochschild homology $HH(O)$ of $O$ is the homology of such complex. 
\end{defi}
\begin{prop} \cite{GV}
The bracket and the star product induce a Gerstenhaber algebra structure 
on the Hochschild homology of an operad with multiplication in vector spaces. 
\end{prop}

Since the operad describing Gerstenhaber algebras is the homology of the little $2$-discs operad $D_2$, Deligne asked his famous question, now known as  the Deligne conjecture, whether the homological action could be induced by an action of (singular) chains of the little discs $C_*(D_2)$ on the Hochschild complex. Many authors 
proved that indeed there was a natural action of a suitable operad quasi-isomorphic to $C_*(D_2)$ on the Hochschild complex.

If we work instead with operads with multiplications
in {\em chain complexes} then the Deligne conjecture holds 
for the {\em normalized} Hochschild complex.
In this context we say that an operad $O$ in chain complexes has a unital multiplication if we have a morphism of operads $Ass_* \to O$ , where 
$Ass_*$ is the operad describing {\em unital} associative algebras. This 
latter operad is also isomorphic as non-symmetric operad to the homology  $H_*(D_1)$ of the little 1-discs. The image of the generator in $Ass_*(0)$
defining the unit is an element $u \in O(0)$.

\begin{defi}
The normalized Hochschild complex of a chain operad with 
(unital) multiplication 
is the subcomplex of the (full) Hochschild complex consisting of those 
elements $x \in O(k),\, k \in \N$ such that $x \circ_i u =0$ for all $1\leq i \leq k$.
\end{defi}

\begin{prop} {\rm (McClure-Smith)} \cite{MS}
The normalized Hochschild complex 
of a chain operad $O$ with unital multiplication
has an action of an operad
quasi-isomorphic to the singular chain operad of the little discs
$C_*(D_2)$.
\end{prop}

It is crucial that the normalized Hochschild complex of a chain operad with unital multiplication 
$O$ can be seen also as (co)normalization of a cosimplicial chain complex
$O^*$ defined from the operad $O$ in a manner completely analogous as in the topological category (section \ref{two}). 
We recall that the (co)normalization of a cosimplicial chain complex $O^*$ is the chain complex
of cosimplicial maps $\Delta^* \ot \Z \to O^*$, with differential induced 
by the cosimplicial chain complex $\Delta^* \ot \Z$.
This construction is the algebraic analog of the totalization of a cosimplicial space.
Thus Theorem \ref{ms}  
can be seen as a topological analog of the Deligne conjecture.
We make this analogy precise in the following statement. 
  
\begin{prop}
Let $O$ be a topological operad with multiplication. The 
Hochschild homology of the operad $C_*(O)$ of singular chains 
on $O$ is isomorphic to the homology of $\,\hotot(O^*)$. The bracket and the star product under the isomorphism 
$HH(C_*(O)) \cong H_*(\hotot(O^*))$
correspond respectively to the Browder operation and the Pontrjagin
product.
\end{prop}

The Gerstenhaber algebra structure interacts well with a spectral sequence computing the homology of $\hotot(O^*)$, the Bousfield
spectral sequence.

\begin{prop} \cite{Bousfield}
Given a cosimplicial space $K^*$, there is a second quadrant spectral sequence
computing the homology of $\hotot{K^*}$. 
Its $E^1$-term  
is $E^1_{-p,q}=H_q(K^p)$,
with the differential  
$\sum_{i=0}^{p+1}(-1)^i d_*^{i}:H_q(K^p) \to H_q(K^{p+1}).$
\end{prop}

The filtration giving the spectral sequence is the 
decreasing filtration by cosimplicial degree in the normalization of
$C_*(K^*)$. 

\begin{prop}
Let $O$ be a topological operad with multiplication. Then the Bousfield spectral sequence for $H_*(\hotot O^*)$ is a spectral sequence of
Gerstenhaber algebras with bracket
$$[\_,\_]: E^r_{-p,q} \ot E^r_{-p',q'} \to E^r_{-p-p'+1,q+q'}$$ and product
$$\_*\_ :E^r_{-p,q} \ot E^r_{-p',q'} \to E^r_{-p-p',q+q'}.$$
The $E_2$-term is the Hochschild homology of the homology operad $H_*(O)$
as a Gerstenhaber algebra.
\end{prop}

\begin{proof}

The star product sums filtration indices on elements in $C_*(O)$.
The bracket $[x,y]$ sits in the  $(m+n-1)$-th filtration term
if $x$ sits in the $m$-th term and $y$ in the $n$-th term.

\end{proof}
The Bousfield spectral sequence does not always converge, but it
does for $K^*=\K_n^*$ or $K^*=\K_n^* \sd S^{n-1}$, as observed by Sinha
\cite{sinha}. 
Arone, Lambrechts and Volic have recently announced a proof
that in these two cases (for $n>3$) the spectral sequence collapses at the $E^2$-term
over the rational numbers \cite{LTV}.
A key ingredient in their proof is a result by Kontsevich showing the formality of the little $n$-discs operad \cite{Kontsevich}, in the sense that
the chain operad $C_*(D_n,\Q)$ is quasi-isomorphic to its homology $H_*(D_n,\Q)$.
The same idea can be used to show
that for $K=\K_n$ there are no extension issues, in the sense that the $E^2$-term is isomorphic to $H_*(Emb'_n,\Q)\cong H_*(\hotot(K_n))$ as 
a Gerstenhaber algebra. We will not need these collapse results here because
in low degree the spectral sequence must collapse and there are no
extension issues.

The $E^2$ term is the Hochschild homology 
of the little $n$-discs operad homology $H_*(D_n)$,
 and has been extensively studied by Turchin \cite{T}.

\

As we have seen the operad $H_*(D_n)$ is generated by a product
$x_1\cdot x_2 \in H_0(D_n(2))$ and a bracket $\{x_1,x_2\} \in H_{n-1}(D_n(2))$.
We use different symbols to avoid confusion with the product and the bracket in the Hochschild complex.

\

{\em Proof of theorem \ref{not2}}.
If $p:Emb'_n \to Emb_n$ is homotopic to a double loop map then it should induce on homology a homomorphism of Gerstenhaber algebras.
We will show that this is not the case because the kernel of $p_*$ is not an ideal with respect to the bracket.

We are considering the case $n$ odd and $n>3$ over rational coefficients.
The lowest dimensional class in the $E^2$-term for $Emb'_n$ is
the element $\alpha=\{x_1,x_2\} \in E_{-2,n-1}$. There is no class that can kill it,  
so this element survives and represents the generator of $H_{n-3}(Emb'_n)\cong \Q$ coming from
the factor $\Omega^2 S^{n-1}$ with respect to the splitting 
$Emb'_n \simeq Emb_n \x \Omega^2 S^{n-1}$.
For similar reasons $H_{2n-6}(Emb'_n) \cong \Q^2$ is generated by
the surviving elements $\beta=\{x_1,x_3\}\cdot\{x_2,x_4\}$ and 
$\alpha^2 = \alpha*\alpha = \{x_1,x_2\}\cdot \{x_3,x_4\}$.

The cosimplicial inclusion $p^*:\K_n^* \to \K_n^* \sd S^{n-1}$ induces 
a morphism of spectral sequences, and on homotopy totalizations
gives a map that we can identify to $p:Emb'_n \to Emb_n$.

The lowest dimensional class in the $E^2$-term for $Emb_n \simeq 
\hotot(\K_n^* \sd S^{n-1})$ is the image $E^2(p)(\beta)$.  
This class
survives to a class $p_*(\beta)$ generating $H_{2n-6}(Emb_n)\cong \Q$.

The computation by Turchin given in
formula 2.9.21 of \cite{T} 
indicates that
the $E^2$-term for $Emb'_n$ in dimension $3n-8$ has two generators,
$[\alpha,\beta]$ and $[\alpha,\alpha^2]=2\alpha[\alpha,\alpha]$, that survive, so that $H_{3n-8}(Emb'_n)\cong \Q^2$.
The $E^2$-term for $Emb_n$ in the same dimension has one generator,
$E^2(p)[\alpha,\beta]$, that survives so that
$H_{3n-8}(Emb_n) \cong \Q$ is generated by $p_*([\alpha,\beta])\neq 0$.
But by dimensional reason $p_*(\alpha)=0$, so the bracket is not preserved
by $p_*$.

Thus $p$ is not a double loop map.
Actually this shows more: there is no double loop space splitting
$Emb'_n \simeq Emb_n \times \Omega^2 S^{n-1}$. 

Now $p$ factors through $fEmb_n$ via a double loop map
 $p':Emb'_n \to fEmb_n$,
that is induced by the operad inclusion $\K_n \to f\K_n$.
This map $p'$ can be identified to the map $Emb_n \times \Omega^2 S^{n-1} \to
Emb_n \times \Omega SO(n-1)$ induced by looping the holonomy
 $\Omega S^{n-1} \to SO(n-1)$. It is well known that $p'_*(\alpha)$ is non-trivial
so by the same reason
the projection $fEmb_n \to Emb_n$ is not
 a double loop map. $\Box$

\medskip

We remark that the obstruction argument does not work rationally for $n$ even because in that
case there is a Gerstenhaber structure on
the $E^2$-term for $Emb_n$  such that $E^2(p)$
 is a map of Gerstenhaber algebras.
Namely additively this $E^2$-term
is identified to
the Hochschild homology of the Batalin-Vilkovisky operad $BV_n$ \cite{T}.
This operad in vector spaces
is the semidirect product of the little $n$-discs homology
$H_*(D_n)$ and the exterior algebra on a generator in dimension $(n-1)$ 
\cite{SW}. 
Then $E^2(p)$ is naturally the
map of Gerstenhaber algebras induced in Hochschild homology by the operad inclusion
$H_*(D_n) \to BV_n$.

\

However only for $n=2$ the operad $BV_2$ is the homology of a topological operad, the framed little $2$-discs operad $fD_2$.
It might be possible that torsion operations like Dyer-Lashof operations still give obstructions
to a double loop structure on the projection $fEmb_n \to Emb_n$ for $n$ even.

\section{String topology of knots} \label{five}

We will show that the suspension spectrum of the space of knots in a sphere,
suitably desuspended, is an $E_2$-ring spectrum, proving theorem \ref{sphere}.

We proved this for $n=3$ in our joint paper with Kate Gruher \cite{GS}.
The original proof was based on the work by Budney, and on a
generalized approach to string topology, 
expanding on fundamental ideas by Chas-Sullivan \cite{CS} and Cohen-Jones \cite{CJ}. 
Now, knowing that $Emb_n$ is a double loop space, we can produce
a proof for $n>3$.
We recall some terminology and we refer to \cite{GS} for details.
We recall that an $E_2$-operad is a topological operad weakly equivalent to the little $2$-discs operad.
Similarly an $E_2$-operad spectrum is an operad in the category of (symmetric) spectra weakly equivalent to
the suspension spectrum of the little $2$-discs operad. For us an $E_2$-ring spectrum will be an algebra over an $E_2$-operad spectrum 
in the weak sense, meaning that the associativity and unit axioms hold up to homotopy.
Given a manifold $M$ with tangent bundle $TM$ we denote by $-TM$ the opposite virtual bundle.

\begin{lem} (Gruher-S.) \label{gs}

Let $X$ be an algebra over an $E_2$-operad $O$, $G$ a compact Lie group and 
$H \subset G$ a closed subgroup. Suppose that $H$ acts on $X$ and the structure maps 
are $H$-equivariant.
Let $p:G \x_H X \to G/H$ be the projection. 
Then the Thom spectrum of the virtual bundle $p^*(-T(G/H))$ over $G \x_H X$ is an 
$E_2$-ring spectrum. 
\end{lem}

Let $Emb(S^1,S^n)$ be the space of smooth embeddings $S^1 \to S^n$.

{\em Proof of thm \ref{sphere}}

It is convenient to use the model for the space of long knots $Emb_n$ given by embeddings 
of the interval into a cylinder
$I \to D_{n-1} \times I$ , with 
$D_{n-1}$ the unit $(n-1)$-disc,  sending 
0 and 1 to $(0,0)$ and $(0,1)$ respectively with tangents directed along the positive direction of the long axis, namely the last coordinate axis.
There is a natural action by $SO(n-1)$ on $Emb_n$ rotating long knots around the long axis.

We have seen in section \ref{two} that $Emb_n$ is weakly equivalent to the homotopy fiber $F$ of $\hotot \K_n^* \to \hotot B_n^*$, by a 
sequence of weak equivalences $$F \to F' \to 
\hotot(\K_n^* \sd S^{n-1}) \leftarrow Emb_n \, ,$$
where
$F'$ is the homotopy fiber of 
$\hotot(K_n^* \sd S^{n-1}) \to \hotot(B_n^* \sd S^{n-1})$.
Actually all maps in the sequence are $SO(n-1)$-equivariant maps between 
$SO(n-1)$-spaces.
Namely the action of $SO(n-1) \subset SO(n)$ on 
$S^{n-1}$ makes $B_n^*$ and $B_n^* \sd S^{n-1}$ into cosimplicial 
$SO(n-1)$-spaces, such that respectively
$\K_n^*$ and $\K_n^* \sd S^{n-1}$ are $SO(n-1)$-invariant cosimplicial subspaces.
Thus the induced maps on homotopy totalizations are $SO(n-1)$-equivariant.
Moreover it is easy to see that the evaluation $Emb_n \to \hotot(K_n^* \sd S^{n-1})$ is $SO(n-1)$-equivariant. Thus $SO(n+1) \x_{SO(n-1)} Emb_n $
is weakly equivalent to $SO(n+1) \x_{SO(n-1)} F$.
As observed by Budney and Cohen \cite{BC}
there is a homotopy equivalence $Emb(S^1,S^n) \simeq SO(n+1) \times_{SO(n-1)} Emb_n$. 
We obtain then a weak equivalence
$Emb(S^1,S^n) \simeq SO(n+1) \times_{SO(n-1)} F $.
 
The $SO(n-1)$-action makes $\K_n$ and $B_n$ into operads
in the category of based $SO(n-1)$-spaces. Thus
the homotopy totalizations
of $\K_n^*$ and $B_n^*$  
are algebras over the operad $\tD_2$ in the category of $SO(n-1)$-spaces, where a trivial $SO(n-1)$-action is assumed on
$\tD_2$.
The inclusion $\hotot \K_n^* \to \hotot B_n^*$ respects this structure,
so that the homotopy fiber $F$ is also an algebra over $\tD_2$
in $SO(n-1)$ spaces.
By lemma \ref{gs},
with $G=SO(n+1),\, H=SO(n-1)$ and $O=\tD_2$,
$(SO(n+1) \times_{SO(n-1)} F)^{-T(SO(n+1)/SO(n-1))}$ is an $E_2$-ring spectrum.
But $SO(n+1)/SO(n-1)$ is (stably) parallelizable
and has dimension $2n-1$,
so that $$(SO(n+1) \times_{SO(n-1)} F)^{-T(SO(n+1)/SO(n-1))} \simeq \Sigma^{1-2n}\Sigma^{\infty} Emb(S^1,S^n)_+$$ is an $E_2$-ring spectrum. $\Box$

The following corollary has been proved independently by
Abbaspour-Chataur-Kallel, who describe also a BV-algebra structure.
\begin{cor}
The homology $H_{*+2n-1}(Emb(S^1,S^n))$ has a natural structure of Gerstenhaber
algebra.
\end{cor}


\begin{thebibliography}{99}
\bibitem{Bu3} R. Budney, The framed little discs operad is cyclic, 
math.AT/0607490
\end{thebibliography}
\end{document}